\def\resume{\if@twocolumn
\section*{R\'esum\'e}
\else \small 
\begin{center}
{\bf R\'esum\'e\vspace{-.5em}\vspace{0pt}} 
\end{center}
\quotation 
\fi}
\def\endresume{\if@twocolumn\else\endquotation\fi}

\documentclass[11pt]{article} 
\usepackage{latexsym}

\newcommand{\qed}{\ \hfill\mbox{$\Box$}\vspace{\baselineskip}}

\newtheorem{theorem}{Theorem}[section]

\newtheorem{proposition}[theorem]{Proposition}
\newtheorem{definition}{Definition}
\newtheorem{corollary}[theorem]{Corollary}

\newenvironment{proof}{\noindent {\bf Proof:}}{{\qed}}

\newcommand{\block}{{\rm\bf B}}
\newcommand{\Nnn}{\mbox{\bf N}}

\newcommand{\iv}{\ensuremath{{\cal I}}}
\newcommand{\lv}{\ensuremath{\ell}}
\newcommand{\LV}{\mbox{$L$}}
\newcommand{\rank}{\rho}
\newcommand{\1}{\hat{1}}
\newcommand{\0}{\hat{0}}

\begin{document}

\title{Generalizations of Eulerian partially ordered sets, flag numbers, and
       the M\"{o}bius function}


\author{Margaret M. Bayer\thanks{This research was supported by University of
        Kansas General Research allocation \#3552.} \\
	Department of Mathematics\\
	University of Kansas\\
	Lawrence KS 66045-2142\\ \and
        G\'abor Hetyei\thanks{On leave from the
	Alfr\'{e}d R\'{e}nyi Institute of Mathematics,
	Hungarian Academy of Sciences.
	Partially supported by Hungarian National Foundation for
	Scientific Research grant no. F 032325.\hfill\break
}\\
	Department of Mathematics\\
        University of North Carolina at Charlotte\\
	Charlotte, NC  28223}

\date{September 2000}

\maketitle

\begin{abstract}
A partially ordered set is {\em $r$-thick} if every nonempty open 
interval contains at least $r$ elements.
This paper studies the flag vectors of graded, $r$-thick posets and shows the
smallest convex cone containing them is isomorphic to the cone of flag
vectors of all graded posets.
It also defines a $k$-analogue of the M\"{o}bius function and $k$-Eulerian
posets, which are $2k$-thick.
Several characterizations of $k$-Eulerian posets are given.
The generalized Dehn-Sommerville equations are proved for flag vectors of
$k$-Eulerian posets.
A new inequality is proved to be valid and sharp for rank 8 Eulerian posets.
\end{abstract}

\begin{resume}
Un ensemble partiellement ordonn\'e est {\em $r$-\'epais} si chacun de
ses intervals ouverts non-vides contient au moins $r$ \'el\'ements. Dans cet
article nous \'etudions les vecteurs $f$ drapeau des ensembles
partiellement ordonn\'es gradu\'es $r$-\'epais. Nous d\'emontrons que le
c\^{o}ne le plus petit contenant ces vecteurs est isomorphe au c\^{o}ne   
des vecteurs $f$ drapeau des ensembles partiellement ordonn\'es
gradu\'es quelconques. Nous d\'efinissons aussi un $k$-analogue de la
fonction de M\"{o}bius et des ensembles partiellement ordonn\'es
$k$-Eul\'eriens qui sont $2k$-\'epais. Nous caract\'erisons les
ensembles partiellement ordonn\'es Eul\'eriens de plusieurs mani\`eres,
et montrons la g\'en\'eralisation des \'equations de Dehn-Sommerville pour
le vecteur $f$ drapeau d'un ensemble partiellement ordonn\'e
$k$-Eul\'erien. 
Nous montrons une nouvelle inegalit\'e optimale pour les ensembles
partiellement ordonn\'es Eul\'eriens de rang~8.     
\end{resume}

\section{Introduction}
In this paper we study certain classes of graded partially ordered sets
(posets), defined by conditions on the sizes of rank sets in intervals.
We are concerned with numerical parameters of the posets, in particular,
flag vectors and the M\"{o}bius function.

A {\em graded poset} $P$ is a finite 
partially ordered set with a unique minimum element
$\0$, a unique maximum element $\1$, and a {\em rank function} $\rank:
P\longrightarrow \Nnn$ satisfying $\rank (\0)=0$, and
$\rank(y)-\rank(x)=1$ whenever $y\in P$ covers $x\in P$.
The {\em rank $\rank(P)$ of a graded poset $P$} is the rank of its maximum 
element. Given a graded poset $P$
of rank $n+1$ and a subset $S$ of $\{1,2,\ldots,n\}$ (which we abbreviate
as $[1,n]$), define the {\em
$S$--rank--selected subposet of $P$} to be the poset
\[P_{S} = \{ x \in P\::\: \rank(x) \in S\} \cup \{ {\0},{\1}\}.\]
Denote by $f_S (P)$ the number of maximal chains of $P_S$. Equivalently,
$f_S(P)$ is the number of chains $x_1<\cdots<x_{|S|}$ in $P$ such that 
$\{\rank(x_1),\ldots,\rank(x_{|S|})\}=S$.
(Call such a chain an {\em $S$-chain} of $P$.)
The vector 
$\left(f_S (P)\::\: S\subseteq [1,n]\right) $
is called the {\em flag
$f$-vector} of $P$. Whenever it does not cause confusion, we write 
$f_{s_1\,\ldots\, s_j}$ rather than $f_{\{s_1,\ldots,s_j\}}$; in particular,
$f_{\{i\}}$ is always denoted $f_i$. 

In the last twenty years there has grown a body of work on numerical conditions
on flag vectors of posets and complexes, especially those arising in geometric
contexts.
A major recent contribution is the determination of the closed
cone of flag vectors of all graded posets by Billera and Hetyei
(\cite{Billera-Hetyei}).
In \cite{Bayer-Hetyei} the authors study the closed cone of flag vectors
of {\em Eulerian} posets.
These are graded  posets for which every (closed) interval has the same number 
of elements of even rank and of odd rank.

A poset is {\em $r$-thick\/} if every nonempty open interval has at least 
$r$ elements.
Thus, every poset is $1$-thick, and Eulerian posets are $2$-thick.  
In the first part of this paper we show that the closed cone of flag vectors
of $r$-thick posets is linearly equivalent to the Billera-Hetyei cone, the
closed cone of flag vectors of all graded posets.

The second part of the paper defines a $k$-analogue of the M\"{o}bius function
and $k$-Eulerian posets (which are $2k$-thick).
We show that the generalized Dehn-Sommerville equations of \cite{Bayer-Billera}
transfer to $k$-Eulerian posets.
These equations have a particularly nice representation in terms of the 
$L^k$-vector, introduced here as a relative of the $cd$-index.
The results of this paper can be used to find inequalities valid for flag
vectors of Eulerian posets.  
In the last section we give as an example a new, sharp inequality for rank~8
Eulerian posets.

\part{$r$-thick posets}

\section{Flag vectors of arbitrary graded posets}

We describe first the cone of flag vectors of all graded posets.
This is due to Billera and Hetyei (\cite{Billera-Hetyei}).

An {\em interval system on $[1,n]$} is any set of subintervals of $[1,n]$
that form an antichain (that is, no interval is contained in another).
A set $S\subseteq [1,n]$ {\em blocks} the interval
system $\iv$ if it has a nonempty intersection with every $I\in \iv$. 
The family of all subsets of $[1,n]$ blocking $\iv$ is denoted by
$\block_{[1,n]}(\iv)$. The main result of \cite{Billera-Hetyei} is the
following.  

\begin{theorem}
\label{T_lgmain}
An expression $\sum _{S\subseteq [1,n]}
a_S f_S(P)$ is nonnegative for all graded posets $P$ of rank $n+1$  if and only
if 
\begin{equation}
\label{E_lgmain}
\sum _{S\in \block_{[1,n]}(\iv)} a_S\geq 0 \quad \mbox{for every interval system $\iv$ on $[1,n]$.}
\end{equation}
\end{theorem}
Here is an outline of the proof from \cite{Billera-Hetyei}.
The proof of the necessity of the condition (\ref{E_lgmain})
involves constructing for every interval system $\iv$ on $\{1,2,\ldots,n
\}$ a family of posets $\{P(n,\iv,N)\::\: N\in \Nnn\}$  of rank $n+1$
such that 
$$\lim _{N\longrightarrow \infty} \frac{1}{
{f_{[1,n]}(P(n,\iv,N))}} \sum _{S\subseteq
[1,n]}a_S f_S(P(n,\iv,N)))=\sum _{{S\in \block_{[1,n]}(\iv)}} a_S.$$

For the other implication,
let $P$ be an arbitrary graded poset, and assume that its
Hasse-diagram is drawn in the plane. Given an interval $[x,y]$ of $P$, let
$\phi(x,y)$ denote the leftmost atom in $[x,y]$. (If $y$ covers $x$ 
then set $\phi(x,y)= y$.) The operation $\phi$ has the
following crucial property: 
\begin{equation}
\label{E_phi}
\mbox{if $p\in [x,y]\subseteq [x,z]$ and $p=\phi([x,z])$ then
$p=\phi([x,y])$.}
\end{equation}  
For every $S\subseteq [1,n]$ and $i\in [1,n]$ define 
$M_S (i)$  to be the smallest $j\in [i,n+1]$ such that $j\in S\cup\{n+1\}$
Consider the set of maximal chains
$$F_S=\left\{\0=p_0<p_1<\cdots<p_n<p_{n+1}=\1 \::\: \forall
i\in[1,n], p_i=\phi ([p_{i-1},p_{M_S(i)}])\right\}.
$$
It is easy to verify that $F_S$ contains exactly $f_S(P)$ elements. 
Moreover, there is a way of associating a family of intervals 
$\iv_C$ to every maximal chain $C=
\{\0=p_0<p_1<\cdots<p_n<p_{n+1}=\1\}$
such that $C$ belongs to $F_S$ if and only if $S$ blocks $\iv_C$. 
The fact that one may find such a family of
intervals is a direct consequence of property (\ref{E_phi}).

\section{Flag vectors of $r$-thick posets}
\label{S_main}

It is easy to expand any graded poset to obtain an $r$-thick poset.
Let $P$ be a graded poset of rank $n+1$.
Write $D^rP$ for the poset obtained from $P$ by replacing every
$x\in P\setminus\{\0,\1\}$ with $r$ elements $x_1$, $x_2$, \ldots $x_r$,
such that $\0$ and $\1$
remain the minimum and maximum elements of the partially ordered set, and
$x_i< y_j$ if and only if $x<y$ in $P$.
The poset $D^rP$ is an $r$-thick graded poset of rank $n+1$.
Clearly $f_S(D^rP)=r^{|S|} f_S(P)$.

\begin{theorem}
\label{T_main} 
For every positive integer $r$,
$\sum_{S\subseteq [1,n]} a_S f_S(P)\geq 0$
for every graded poset $P$ of rank $n+1$ if and only if 
$\sum_{S\subseteq [1,n]} a_S r ^{n-|S|}f_S (Q)\geq 0$ for
every $r$-thick poset $Q$ of rank $n+1$. 
\end{theorem}

\begin{proof}
First assume $\sum_{S\subseteq [1,n]} a_S r ^{n-|S|}f_S (Q)\geq 0$ for
every $r$-thick poset $Q$ of rank $n+1$. 
Let $P$ be any graded poset of rank $n+1$.
Since $D^rP$ is $r$-thick, 
\begin{eqnarray*}
0&\leq &\sum_{S\subseteq [1,n]} a_S r^{n-|S|}f_S (D^rP)\\
 &=& \sum_{S\subseteq [1,n]} a_S r^{n-|S|}r^{|S|}f_S (P)\\
 &=& \sum_{S\subseteq [1,n]} a_S r^n f_S (P).
\end{eqnarray*}
Dividing by $r^n$ gives the desired inequality for all graded posets.

Now assume $\sum_{S\subseteq [1,n]} a_S f_S(P)\geq 0$
for every graded poset $P$ of rank $n+1$.
Let $Q$ be an $r$-thick poset of rank $n+1$.
For each rank $i$, fix a total order of the elements of $Q$ of rank $i$.
Given an interval $[x,y]$ of $Q$ of rank at least 2, let
$\phi(x,y)$ denote the set of the first $r$ atoms in $[x,y]$. 
(If $y$ covers $x$, set $\phi(x,y) = \{y\}$.)
 
The operation $\phi$ satisfies the following:
\begin{equation}
\label{E_phi2}
\mbox{if $p\in [x,y]\subseteq [x,z]$ and $p\in\phi([x,z])$ then
$p\in\phi([x,y])$.}
\end{equation}  
Let
$$F_S=\left\{\0=p_0<p_1<\cdots<p_n<p_{n+1}=\1 \::\: \forall
i\in[1,n], p_i\in\phi ([p_{i-1},p_{M_S(i)}])\right\}.
$$ 
How many sequences are in the set $F_S$?
Given any $S$-chain of $Q$, extend it to sequences in $F_S$ one rank
at a time.
Having fixed $p_0$ through $p_{i-1}$ ($1\le i\le n$), if $i\not\in S$,
then there are exactly $r$ choices for $p_i$.
Thus $|F_S|=r^{n-|S|}f_S(Q)$.

To each maximal chain $C$: $\0=p_0<p_1<\cdots<p_n<p_{n+1}=\1$ of $Q$ is assigned
an interval system as follows. 
For $1\le i\le n$, let $\psi(C,i)$ be the largest $j$ such that 
$p_i\in\phi(p_{i-1},p_j)$.
Let $\iv'_C=\{[i,\psi(C,i)]:\, \mbox{$1\le i\le n$, $\psi(C,i)\ne n+1$}\}$, and 
let $\iv_C$ be the interval
system consisting of minimal intervals in $\iv'_C$.
We show $C$ belongs to $F_S$ if and only if $S$ blocks $\iv_C$.
Suppose $ C$: $\0=p_0<p_1<\cdots<p_n<p_{n+1}=\1$ is in $F_S$.
Then for all $i$, $p_i\in \phi([p_{i-1},p_{M_S(i)}])$, so by the maximality of 
$\psi(C,i)$, $\psi(C,i)\ge M_S(i)$.
So for all $i$ the interval $[i,\psi(C,i)]$ contains the element $M_{S(i)}$ of 
$S$.
Thus $S$ blocks $\iv_C$.
For the reverse implication, suppose $C$ is a maximal chain of $Q$ and $S$
blocks $\iv_C$.
Let $1\le i\le n$ and $[i,\psi(C,i)]\in \iv_C$.
Since $S$ blocks $\iv_C$, $S\cap [i,\psi(C,i)]$ contains an element $s$.
So $M_{S(i)}\le s \le \psi(C,i)$.
Apply condition (\ref{E_phi2}): $p_i\in[p_{i-1},p_{M_S(i)}]\subseteq
[p_{i-1},p_{\psi(C,i)}]$ and $p_i\in\phi([p_{i-1},p_{\psi(C,i)}])$, so
$p_i\in\phi([p_{i-1},p_{M_{S(i)}}])$.
Thus $C$ is in $F_S$.

Given a system of intervals $\iv$ denote by $f_{\iv}$ the number
of those maximal chains $C$ of $Q$ for which 
$\iv_C=\iv$. 
(Note that $f_\iv$ depends not only on 
$Q$ but also on the ordering of the elements of each rank.)
Then
\begin{eqnarray*}
\sum_{S\subseteq [1,n]} a_S r^{n-|S|} f_S (Q)
&=&\sum_{S\subseteq [1,n]} a_S |F_S|
=\sum_{S\subseteq [1,n]} a_S \sum _{S\in\block_{[1,n]} (\iv)} f_{\iv}\\
&=&\sum _{\iv} f_{\iv} \sum _{S\in\block_{[1,n]} (\iv)} a_S.\\ 
\end{eqnarray*}
By Theorem \ref{T_lgmain} the sums $\sum _{S\in\block_{[1,n]} (\iv)} a_S$
are all nonnegative, and so 
$$\sum_{S\subseteq [1,n]} a_S r^{n-|S|} f_S (Q)\geq 0.$$ 

\vspace*{-30pt}

\end{proof}

Let ${\cal C}_{r,n+1}$ be the smallest closed convex cone containing the flag
vectors of all $r$-thick posets of rank $n+1$.
\begin{corollary}
\label{cone-cor}
For all positive integers $q$ and $r$, the invertible linear transformation
$\alpha_{q,r}: {\bf Q}^{2^n}\rightarrow {\bf Q}^{2^n}$ defined by
$\alpha_{q,r}((x_S))=((r/q)^{|S|}x_S)$ maps ${\cal C}_{q,n+1}$ onto 
${\cal C}_{r,n+1}$.
\end{corollary}

To determine if a graded poset is $r$-thick, it is enough to check that
between every $x$ and $y$ with $x < y$ and $\rank(y)-\rank(x)=2$, there
are at least $r$ elements.
The definition of $r$-thick posets can then be generalized by allowing
the lower bound $r$ to vary through the levels of the poset.
The results of this section have straightforward analogues in that
context.

\part{$k$-Eulerian posets}

\section{The $k$-M\"obius function}

\begin{definition}
{\em The {\em M\"obius function} of a graded poset
$P$ is defined recursively for any subinterval of $P$ by the formula

\[\mu ([x,y])=                                                                  
\left\{                                                                         
\begin{array}{cl}                                                               
1 & \mbox{if $x=y$},\\                                                          
-\sum_{x\le z<y} \mu([x,z]) &\mbox{otherwise}.\\   
\end{array}                                                                     
\right.                                                                         
\] 

A graded poset $P$ is {\em Eulerian} if the M\"obius function of every interval 
$[x,y]$ is given by $\mu([x,y]) = (-1)^{\rank(x,y)}$. }
\end{definition}

(Here $\rank(x,y)=\rank([x,y])=\rank(y) - \rank(x)$.)                        


See \cite{Stanley-Eulerian} for a survey of Eulerian posets.
The first characterization of all linear equalities holding 
for the flag vectors of all Eulerian posets was given by Bayer and 
Billera in \cite{Bayer-Billera}. 
\begin{theorem}[Bayer and Billera]
\label{T_BB}
For every Eulerian poset of rank \mbox{$n+1$}, every subset
$S\subseteq [1,n]$, and every maximal interval $[i,\ell]$  
of $[1,n]\setminus S$,
\[\left((-1)^{i-1}+(-1)^{\ell+1}\right)f_S(P)+\sum _{j=i}^\ell (-1)^j
f_{S\cup\{j\}}(P)=0.\]
Furthermore, every linear equality holding for the flag vector of all Eulerian 
posets of rank $n+1$ is a consequence of these equations.
\end{theorem}

Next we present generalizations of the M\"obius function and of Eulerian
posets.

\begin{definition}
{\em
The {\em $k$-M\"obius function} of a graded poset is defined 
recursively by the formula 
$$\mu_k ([x,y])= 
\left\{
\begin{array}{cl}
1 & \mbox{if $x=y$},\\
-1-\frac{1}{k}\sum_{x<z<y} \mu_k([x,z]) &\mbox{otherwise}.\\
\end{array}
\right.
$$}
\end{definition}

The following proposition gives the $k$-M\"{o}bius function of a poset $P$
as a $k$-analogue of the reduced Euler characteristic of the order complex 
of $P$.
It is a generalization of Philip Hall's theorem, and 
is easy to prove by induction.
\begin{proposition}
\label{E_MqEc}
If $P$ is a graded poset of rank $n+1$, then
$$ \mu_k (P)=-\sum_{S\subseteq [1,n]} (-\frac{1}{k})^{|S|} f^{n+1}_S(P).  $$
\end{proposition}

A graded poset is {\em $k$-Eulerian} if for every interval 
$[x,y]\subseteq P$, $\mu_k ([x,y])=(-1)^{\rank(x,y)}$. 
Note that $1$-Eulerian is the same as Eulerian.
The following proposition follows easily from the definitions.
\begin{proposition}
\label{Eulfvec}
If $P$ is a $k$-Eulerian poset of rank $n+1$, then
\begin{enumerate}
\item every interval of $P$ is $k$-Eulerian
\item $\displaystyle \sum_{i=1}^n (-1)^{i-1}f_i(P)=k(1-(-1)^n)$
\label{Eulfvec2}
\end{enumerate}
\end{proposition}

The thickening operation introduced in Section~\ref{S_main} connects the
$k$-M\"obius function for different values of $k$.
\begin{proposition}
\label{P_kreplica}
Let $[x,y]$ be an interval of a graded poset $P$ and $\ell$ a positive integer. 
Consider an interval $[x_i,y_j]\subseteq D^\ell  P $ corresponding 
to $[x,y]\subseteq P$. Then 
$$
\mu_k([x,y])=
\mu_{k\ell} ([x_i,y_j]).
$$
\end{proposition}
\begin{proof}
Recall that $f_S(D^\ell P )=\ell^{|S|} f_S(P)$.
Since the interval $[x_i,y_j]$ of  $D^\ell P$ is isomorphic to 
$D^\ell  [x,y] $, the result is obtained by substitution in 
Proposition~\ref{E_MqEc}. 
\end{proof}

\begin{corollary}
\label{q-kq}
A poset $P$ is $k$-Eulerian if and only if $D^\ell P$ is $k\ell$-Eulerian.
\end{corollary}
In \cite{Bayer-Hetyei} a half-Eulerian poset was defined to be a poset $P$
for which $D^2P$ is Eulerian.

Using Proposition \ref{P_kreplica} we can determine exactly the set of those 
$k$'s for which $k$-Eulerian posets exist. 
\begin{theorem}
For every positive integer $n$,
there exists a $k$-Eulerian poset of rank $n+1$ if and only if $k=j/2$ for some 
positive integer $j$.
Moreover, every $k$-Eulerian poset is $2k$-thick.
\end{theorem}
\begin{proof}
The chain $C$ of rank $n+1$ is half-Eulerian.
For every positive integer $j$, $D^jC$ is a $j/2$-Eulerian poset.
On the other hand, by the definition of the function $\mu_k$, 
for an interval $[x,y]$ of rank $2$ in a $k$-Eulerian poset,
$$(-1)^2=\mu_k ([x,y])
=-1-{1/k}\sum_{x<z<y}\mu_{1/k} ([x,z])
=-1-{1\over k}\sum_{x<z<y} (-1).
$$
Therefore $2k$ is the number of elements $z$ strictly between $x$ and $y$. 
Thus, if $P$ is a $k$-Eulerian poset, then $2k$ is a positive integer, and
$P$ is $2k$-thick.
\end{proof}

It is easy to check by induction that a graded poset is half-Eulerian if 
and only if (1) in every interval $[x,y]$ with $\rho(x,y)$ odd,
the number of elements of even rank equals the number of elements of odd rank;
and (2) in every interval $[x,y]$ with $\rho(x,y)$ even, the number of
elements $z$ with $\rho(x,z)$ even is one more than the number of
elements $z$ with $\rho(x,z)$ odd.
This characterization can be used to check that the following 
``vertical doubling'' of an arbitrary graded poset produces a
half-Eulerian poset. Let $P$ be any graded poset with relation $\prec_P$.
Form the set $Q=\{\hat{0},\hat{1}\}\cup\{x_1,x_2 : x\in P\setminus
\{\hat{0},\hat{1}\}\}$.
Define a relation $\prec_Q$ on $Q$ by $u\prec_Q v$ if and only if 
one of the following holds:

$\bullet$ $u=\hat{0}$, $v\in P\setminus\{\hat{0}\}$

$\bullet$ $v=\hat{1}$,  $u\in P\setminus\{\hat{1}\}$

$\bullet$ $u=x_1$ and $v=x_2$ for some $x\in P\setminus\{\hat{0},\hat{1}\}$

$\bullet$ $u=x_i$ and $v=y_j$ for some $x, y\in P\setminus\{\hat{0},\hat{1}\}$,
      with $x\prec_P y$.

\noindent
If $P$ is a rank $n+1$ graded poset, then the resulting poset $Q$ is a rank 
$2n+1$ half-Eulerian poset.

For larger $k$, 
not all $k$-Eulerian posets are obtained by the thickening operation.
For an example, consider the poset $P$ of rank $n+1\ge 3$ having 
elements $x_1$, $x_2$, \ldots, $x_m$ of rank~1, elements 
$y_1$, $y_2$, \ldots, $y_m$ of rank~2, with $x_i < y_j$ if and only if
$i=j$, and one element of each other rank.
It is easy to check that $P$ is half-Eulerian, and so $D^{2k}P$ is
$k$-Eulerian.  
In the Hasse diagram of $D^{2k}P$, the subgraph induced by 
the elements of ranks~1 and~2 consists of $m$ copies of the complete
bipartite graph $K_{2k,2k}$.
Replace this subgraph by any other $2k$-regular bipartite graph on these
elements.
The resulting graph is the Hasse diagram of another $k$-Eulerian poset.
(Note that the only relations changed in the poset are those between
rank~1 and rank~2 elements.)

The definition of $k$-Eulerian, like that of $r$-thick, can be 
generalized by varying the multiplier $k$ with the rank of the elements.
The results of this and the next section can easily be adapted for such
posets.

\section{The flag $L^k$-vector}

A certain transformation of the flag $f$-vector was useful in 
\cite{Stanley-flag}, \cite{Billera-Hetyei},
and \cite{Bayer-Hetyei}. 
It has a natural adaptation to the $k$-Eulerian setting.
\begin{definition}
\label{D_LVq}
{\em
The {\em flag $L^k$-vector} of a graded partially ordered 
set $P$ of rank $n+1$ is the vector 
$(L^{k,n+1}_S(P) \::\: S\subseteq [1,n])$, where
$$L^{k,n+1}_S (P) = (-1)^{n-|S|}
\sum_{T\supseteq [1,n]\setminus S} \left(-{1\over 2k}\right)^{|T|}
f^{n+1}_T(P).$$   
}
\end{definition}
For $k=1/2$ this is the \lv-vector of \cite{Billera-Hetyei};
for $k=1$ this is the ``$ce$-index'' of \cite{Stanley-flag} and the
\LV-vector of \cite{Bayer-Hetyei}.
The formula inverts to give
\begin{equation}
\label{E_fLq}
f^{n+1}_S(P)=(2k)^{|S|}\sum _{T\subseteq [1,n]\setminus S} L^{k,n+1}_T(P).
\end{equation}
The $L^k$-vector ignores the effect of the operator $D^\ell$.
If $P$ is a graded poset of rank $n+1$, then
\begin{equation}
\label{Lq-Lkq}
L^{k\ell,n+1}_S (D^\ell P)=L^{k,n+1}_S(P).
\end{equation}

A set $S\subseteq[1,n]$ is {\em even} if $S$ is a disjoint union of 
intervals of even cardinality.
The parameters $L^{k,n+1}_S$ for even sets $S$ play a special role for 
$k$-Eulerian posets.
The $k$-analogue of Theorem \ref{T_BB} is the following. 
\begin{theorem}
\label{T_BBq}
For every $k$-Eulerian poset $P$ of rank $n+1$, every subset
$S\subseteq [1,n]$, and every maximal interval $[i,\ell]$  
of $[1,n]\setminus S$,
\[k\left((-1)^{i-1}+(-1)^{\ell+1}\right)f_S(P)+\sum _{j=i}^\ell (-1)^j
f_{S\cup\{j\}}(P)=0.\]
Every linear equality holding for the flag vector of all 
$k$-Eulerian posets of rank $n+1$ is a consequence of these equations.

In $L^k$-vector form, these equations are equivalent to the set of 
equations $L^{k,n+1}_S(P)=0$ for all subsets $S\subseteq[1,n]$ that are
not even.
\end{theorem}
Call these equations the {\em generalized Dehn-Sommerville equations},
and denote by $DS_{k,n+1}$ the resulting subspace of ${\bf R}^{2^n}$.

\vspace*{6pt}

\begin{proof}
The fact that the equations (in flag $f$-vector form) hold for all
$k$-Eulerian posets follows from Proposition~\ref{Eulfvec}.
Fix a set $S$ with gap $[i,\ell]$. 
For each $S$-chain identify the 
rank $i-1$ element $x$ and rank $\ell+1$ element $y$, and apply 
equation~(\ref{Eulfvec2}) to the interval $[x,y]$.  
Sum the resulting equations for all the $S$-chains.

Convert the flag $f$-vector equations using equation~(\ref{E_fLq}).
Writing $V=[1,n]\setminus S$ and dividing by $2^{|S|}k^{|S|+1}$, the result is
\begin{equation}
\label{E_BBlq}
\left((-1)^{i-1}+(-1)^{\ell+1}\right) \sum _{T\subseteq V} L^{k,n+1}_T 
+2 \sum _{j=i}^\ell (-1)^j \sum _{T\subseteq V\setminus \{j\}}
L^{k,n+1}_T =0.
\end{equation}
  From this we prove by induction that $L^{k,n+1}_V(P)=0$ (abbreviated as
$L_V=0$) for all noneven sets $V$.
Let $V\subseteq[1,n]$ be any noneven set, and let $[i,\ell]$ be an odd maximal
interval of $V$.
Equation~(\ref{E_BBlq}) gives 
\begin{equation}
\label{L-eqn}
\sum_{T\subseteq V} L_T + 
\sum_{j=i}^\ell (-1)^{j-i+1}\sum_{T\subseteq V\setminus\{j\}}L_T = 0. 
\end{equation}
If $T$ is a noneven proper subset of $V$, then by the induction assumption,
$L_T=0$.
So consider an even subset $T\subseteq V$.
Since the maximal intervals of $T$ contained in $[i,\ell]$ are even,
$[i,\ell]\setminus T = \{j_1,j_2,\ldots,j_t\}$, where $t$ is odd, $j_1-i$ is even,
and, for $2\le p\le t$, $j_p - j_{p-1}$ is odd.
Thus, for $1\le p\le t$, $j_p - i + 1$ has the same parity as $p$.
The coefficient of $L_T$ in (\ref{L-eqn}) is $1+\sum_{p=1}^t (-1)^{j_p-i+1}
=1+\sum_{p=1}^t (-1)^p =0$.
So equation~(\ref{L-eqn}) reduces to $L_V=0$.

Conversely, suppose $L_V=0$ for all noneven sets $V\subseteq[1,n]$.
We show that the equations in (\ref{E_BBlq}) hold.
Let $V\subseteq[1,n]$ and $[i,\ell]$ a maximal interval of $V$.
For $\ell-i$ even, we need to prove equation (\ref{L-eqn}).
(The case of $\ell-i$ odd is similar, and is omitted.)
It suffices to consider the terms $L_T$ with $T$ an even set.
For such $T$, $[i,\ell]\setminus T= \{j_1,j_2,\ldots,j_t\}$ as above, with $t$ odd,
and $j_p-i+1\equiv p \pmod{2}$.
So the coefficient in (\ref{L-eqn}) of $L_T$ is 
\mbox{$1+\sum_{p=1}^t (-1)^{j_p-i+1} =1+\sum_{p=1}^t (-1)^p =0$}.
Thus equation~(\ref{L-eqn}) holds.

To complete the proof, it suffices to show that the linear span of the 
$L^k$-vectors of $k$-Eulerian posets of rank $n+1$ is the subspace of
${\bf R}^{2^n}$ determined by the equations 
$L^{k,n+1}_S(P)=0$ for all subsets $S\subseteq[1,n]$ that are not even.
This can be accomplished by finding a set of linearly independent vectors
in the span of the $L^k$-vectors of $k$-Eulerian posets, one vector
for each even subset $S\subseteq[1,n]$.
In \cite{Billera-Hetyei} Billera and Hetyei constructed, for each interval
system $\cal I$, a sequence of graded posets $P(n,{\cal I}, N)$.
The construction starts with a rank $n+1$ chain, and replicates
intervals of ranks in the poset.
For an even set $S$, let ${\cal I}[S]$ be the set of maximal intervals in $S$.
(For example, for $S=\{1,3,4,7,8,9,10\}$, ${\cal I}[S]=\{ [1], [3,4],[7,10]\}$.)
If $S$ is an even subset of $[1,n]$,
then $P(n,{\cal I}[S], N)$ is half-Eulerian for all $N$.
Furthermore, the sequence of $L^{1/2}$-vectors of these posets satisfies
the following.
Here $m$ is the number of intervals in ${\cal I}[S]$.
$$\lim_{N\longrightarrow \infty} {1\over N^m} L^{1/2,n+1}_T 
\left(P(n,\iv[S], N)\right) 
=\left\{
\begin{array}{ll}
(-1)^{j} &\mbox{if $T$ is the union of $j$ intervals of $S$,}\\
0 & \mbox{otherwise}.\\
\end{array}
\right.
$$
(See \cite{Bayer-Hetyei} for details.)
Using (\ref{Lq-Lkq}), we get for any positive integer $2k$,
$$\lim_{N\longrightarrow \infty} {1\over N^m} L^{k,n+1}_T 
\left(D^{2k}P(n,\iv[S], N)\right) 
=\left\{
\begin{array}{ll}
(-1)^{j} &\mbox{if $T$ is the union of $j$ intervals of $S$,}\\
0 & \mbox{otherwise}.\\
\end{array}
\right.
$$
For fixed $k$ the limiting $L^k$-vectors for each even interval system 
${\cal I}[S]$ are linearly independent, since for each even set $S$, the vector
formed from the sequence
$(P(n,{\cal I}[S], N))$ has $T$-entry 0 for all $T$ not containing $S$.
\end{proof}

A flag vector can by chance lie in the subspace $DS_{k,n+1}$
without the poset being $k$-Eulerian.
However, $k$-Eulerian posets are characterized by the equations holding
locally.
The $k=1$ case of this is in \cite{Bayer-Hetyei}.
The proof requires the convolution of flag operators, defined by Kalai
\cite{Kalai} (see also \cite{Billera-Liu}).
It is defined for the flag numbers by
$f^m_S * f^n_T= f^{m+n}_{S\cup \{m\}\cup (T+m)}$, and is extended by bilinearity
to linear combinations.
For $p^{m+1}$ and $q^{n+1}$ linear combinations of chain operators in ranks
$m+1$ and $n+1$, respectively, their convolution on a rank $m+n$ poset $P$
satisfies
$$p^{m+1} * q^{n+1}(P) = \sum_{x\in P\atop \rank(x)=m}
p^{m+1}([\0,x]) q^{n+1}([x,\1]).$$
Convolution behaves nicely on the flag $L^k$-vector.   
For a rank $m+n+2$ poset $P$,
\begin{equation}
\label{convol}
L^{k,m+1}_S * L^{k,n+1}_T(P)= \sum_{x\in P\atop \rank(x)=m+1}
L^{k,m+1}_S ([\0,x])L^{k,n+1}_T ([x,\1])= 2kL^{k,m+n+2}_{S\cup (T+m+1)}(P).
\end{equation}

\begin{theorem}
\label{P_qEulerLV}
A graded partially ordered set $P$ is $k$-Eulerian if and only if for every 
interval $[x,y]\subseteq P$  of positive even rank
$L^{k,\rank(x,y)}_{[1,\rank(x,y)-1]} ([x,y])=0$.
\end{theorem}

\begin{proof}
Since every interval of a $k$-Eulerian partially ordered set is $k$-Eulerian,
Theorem~\ref{T_BBq} gives that
$L^{k,\rank(x,y)}_{[1,\rank(x,y)-1]} ([x,y])=0$
for all intervals $[x,y]$ of positive even rank.

Now assume that for every interval $[x,y]$ of positive even rank, \linebreak
$L^{k,\rank(x,y)}_{[1,\rank(x,y)-1]} ([x,y])=0$.
Then by equation~(\ref{convol}), for every interval $[x,y]\subseteq P$ and
for every
$S\subseteq [1,\rank(x,y)-1]$ that is not even,
$L^{k,\rank ([x,y])}_S ([x,y])=0$.

For $P$ of rank $n+1$, by Proposition~\ref{E_MqEc},
\begin{eqnarray*}
\mu_k(P)&=&-\sum_{S\subseteq [1,n]} (-\frac{1}{k})^{|S|} f^{n+1}_S(P)
=-\sum_{S\subseteq [1,n]} (-\frac{1}{k})^{|S|} (2k)^{|S|} 
\sum_{T\subseteq [1,n]\setminus S} L^{k,n+1}_T (P)\\
&=&-\sum_{T\subseteq [1,n]} L^{k,n+1}_T (P)
     \sum_{S\subseteq [1,n]\setminus T} (-2)^{|S|}
=-\sum_{T\subseteq [1,n]} \LV^{n+1}_T (P) (-1)^{n-|T|}.
\end{eqnarray*}
Since $L^{k,n+1}_T (P)$ is nonzero only if T is an even set, and then $|T|$ is
an even number,
$$\mu_k(P)=(-1)^{n+1} \sum_{T\subseteq [1,n]} \LV^{n+1}_T (P)
=(-1)^{n+1} f^{n+1}_{\emptyset}(P)=(-1)^{n+1}.$$
The same argument can be repeated for every interval of $P$, showing that
it is a $k$-Eulerian poset.
\end{proof}

Using this result, we get the following curious characterization via the
M\"{o}bius function.
\begin{theorem}
A graded poset $P$ is $k$-Eulerian if and only if the $2k$-M\"obius function of
every interval $[x,y]\subseteq P$ of even rank is zero.
\end{theorem}
\begin{proof}
Let $P$ be a graded poset.
By Corollary~\ref{q-kq} $P$ is Eulerian if and only if $D^2P$ is
$2k$-Eulerian 
if and only if for every interval $[x_i,y_j]$ of $D^2P$ with $\rank(x_i,y_j)
$ even
$$L^{2k,\rank(x_i,y_j)}_{[1,\rank(x_i,y_j)-1]} ([x_i,y_j])=0$$
if and only if for every interval $[x_i,y_j]$ of $D^2P$ with $\rank(x_i,y_j)$
even
$$\sum_{T\subseteq[1,\rank(x_i,y_j)-1]}(-\frac{1}{4k})^{|T|}f_T([x_i,y_j])=0$$
if and only if for every interval $[x,y]$ of $P$ with $\rank(x,y)$ even
$$\sum_{T\subseteq[1,\rank(x,y)-1]}(-\frac{1}{4k})^{|T|}2^{|T|}f_T([x,y])=0$$
if and only if for every interval $[x,y]$ of $P$ with $\rank(x,y)$ even
$$\sum_{T\subseteq[1,\rank(x,y)-1]}(-\frac{1}{2k})^{|T|}f_T([x,y])=0$$
if and only if for every interval $[x,y]$ of $P$ with $\rank(x,y)$ even
$\mu_{2k}([x,y])=0$.
\end{proof}

In particular, a graded poset $P$ is half-Eulerian if and only if 
the (usual) M\"obius function of $[x,y]$ is zero 
for every $[x,y]\subseteq P$ of even rank.

The $L^1$-vector of a graded poset is the vector of coefficients of the
$ce$-index, introduced in \cite{Stanley-flag} as a variation of the
$cd$-index of an Eulerian poset.
(The $cd$-index of an Eulerian poset, due to Fine (see \cite{Bayer-Klapper}), 
is a vector linearly equivalent to the flag vector; it embodies the 
generalized Dehn-Sommerville equations of Theorem~\ref{T_BB}.)
In \cite{Stanley-flag}, Stanley observed that the existence of the $cd$-index 
for a
graded poset is equivalent to the vanishing of the coefficients of $ce$-words
containing an odd string of $e$'s; in our notation this says $L^{1,n+1}_S(P)=0$
for all subsets $S\subseteq[1,n]$ that are not even.
Thus the last part of Theorem~\ref{T_BBq} (as well as the first part) 
is already known for $k=1$. 

The $L^k$-vector for general $k$ can be presented in the same way.
For $P$ any graded poset of rank $n+1$, write a generating function for the
flag $f$-vector as follows:
$$\Upsilon(a,b)=\sum_{S\subseteq[1,n]} f_S u_S,$$ where 
$u_S=u_1u_2\ldots u_n$ with $u_i=a$ if $i\not\in S$ and $u_i=b$ if $i\in S$.
then $$\Upsilon(e,\frac{c-e}{2k})=\sum_{T\subseteq[1,n]} L_T v_T,$$
where
$v_T=v_1v_2\ldots v_n$ with $v_i=c$ if $i\not\in T$ and $v_i=e$ if $i\in T$.
The equations of Theorem~\ref{T_BBq} for $k$-Eulerian posets can then be
rephrased as saying that $\Upsilon(e,(c-e)/(2k))$ is a polynomial in
(the noncommuting expressions) $c$ and $ee$.

\section{The cone of $k$-Eulerian flag vectors}
Theorem~\ref{T_main}, along with the description of the cone of flag 
vectors of general graded posets (\cite{Billera-Hetyei}), can be used to 
generate all the inequalities valid for all $r$-thick posets.
The inequalities for $2k$-thick posets are, in particular, valid for all
$k$-Eulerian posets, but they may not be sharp.
We would like to know the essential inequalities, that is, 
the closed cones of flag vectors of Eulerian and of half-Eulerian
posets.
In \cite{Bayer-Hetyei} these cones are studied and are
completely determined up through rank 7.
(See also \cite{homepage} for data on the cone.)
In the context of this paper, the results can be stated as follows.
\begin{theorem}[\cite{Bayer-Hetyei}]
For rank $n+1\le 7$,
\begin{enumerate}
\item
\label{part1}
the closed cone of flag vectors of half-Eulerian posets of rank $n+1$ is 
the intersection of the cone ${\cal C}_{1,n+1}$ of flag vectors of all
graded posets of rank $n+1$ with the subspace $DS_{1/2,n+1}$ determined by the
half-Eulerian equations of Theorem~\ref{T_BBq};
\item
\label{part2}
the closed cone of flag vectors of Eulerian posets of rank $n+1$ is 
the intersection of the cone ${\cal C}_{2,n+1}$ of flag vectors of all
$2$-thick
graded posets of rank $n+1$ with the generalized Dehn-Sommerville
subspace $DS_{1,n+1}$; and
\item
\label{part3}
the two cones are isomorphic.
\end{enumerate}
\end{theorem}

We do not know if this theorem extends to higher ranks.
However, for all ranks, part~\ref{part1} of the theorem implies 
parts~\ref{part2} and~\ref{part3}.
\begin{theorem}
Let ${\rm CONE}_{k,n+1}$ be the statement,
\begin{quote}
The closed cone of flag vectors of $k$-Eulerian posets of rank~\mbox{$n+1$} is the
intersection of the cone ${\cal C}_{2k,n+1}$ of flag vectors of all
$2k$-thick graded posets of rank $n+1$ with the generalized Dehn-Sommerville
space $DS_{k,n+1}$.
\end{quote}
For all $k\ge 1$ (with $2k$ an integer) and all positive integers $n$,
$$ {\rm CONE}_{1/2,n+1} \Longrightarrow {\rm CONE}_{k,n+1}.$$
\end{theorem}

\begin{proof}
Recall the map $\alpha_{1,2k}$ of Corollary~\ref{cone-cor}; it maps
${\cal C}_{1,n+1}$ onto ${\cal C}_{2k,n+1}$.
Clearly it also maps $DS_{1/2,n+1}$ onto $DS_{k,n+1}$.
So $\alpha_{1,2k}({\cal C}_{1,n+1}\cap DS_{1/2,n+1})=
{\cal C}_{2k,n+1}\cap DS_{k,n+1}$, which contains the cone of $k$-Eulerian
flag vectors.
On the other hand, for any half-Eulerian poset $P$, 
$\alpha_{1,2k}((f_S(P)))=(f_S(D^{2k}P))$, the flag vector of the
$k$-Eulerian poset $D^{2k}P$.
If $\mbox{CONE}_{1/2,n+1}$ holds, then 
${\cal C}_{1,n+1}\cap DS_{1/2,n+1}$ is the cone of half-Eulerian flag vectors,
and its image is contained in the cone of $k$-Eulerian flag vectors.
Thus, if $\mbox{CONE}_{1/2,n+1}$ holds, then 
${\cal C}_{2k,n+1}\cap DS_{k,n+1}$ is exactly the closed cone of flag
vectors of $k$-Eulerian posets.
\end{proof}

Another question raised in \cite{Bayer-Hetyei} on the structure of these 
cones can be answered.
For rank at most 7, all facet inequalities of the half-Eulerian
(and with slight modification, Eulerian)
cone are generated from two basic types of inequalities.
\begin{theorem}[\cite{Bayer-Hetyei}]
\label{ineqs}
Let $S$ and $T$ be disjoint subsets of $[1,n]$, such that  every maximal
interval of the complement of $S$ contains at most one element of $T$.
Then for every rank $n+1$ half-Eulerian poset $P$,
$$\sum_{R\subseteq T} (-1)^{|T\setminus R|}f_{S\cup R}(P)\ge 0.$$
Let $1\le i< j<\ell\le n$.
Then for every rank $n+1$ half-Eulerian poset $P$,
$$f_{i\ell}(P)-f_i(P)-f_\ell(P)+f_j(P)\ge 0.$$
\end{theorem}
Other valid inequalities are obtained by the convolution of inequalities
of these types.
The question arose whether these generate all inequalities valid for the
flag vectors of all half-Eulerian posets.
They do not.
\begin{proposition}
For all half-Eulerian posets $P$ of rank~8,
$$f^8_{1356}(P)-f^8_{135}(P)-f^8_{356}(P)+f^8_{15}(P)-f^8_{16}(P)+f^8_{35}(P)+
f^8_{36}(P)-f^8_3(P)\ge 0,$$
or, in $L^{1/2}$-vector form,
\begin{eqnarray}
\label{rank8ineq}
L^{1/2,8}_{45}(P)+ L^{1/2,8}_{2345}(P)+ L^{1/2,8}_{56}(P)+ L^{1/2,8}_{1256}(P)
-L^{1/2,8}_{2367}(P)   
\nonumber \\ 
{}  -  L^{1/2,8}_{3467}(P)+ L^{1/2,8}_{4567}(P)+
L^{1/2,8}_{124567}(P)&\le & 0.
\end{eqnarray}
This inequality determines a facet of the closed cone of flag vectors
of half-Eulerian posets, and does not follow from the inequalities 
of Proposition~\ref{ineqs}.
\end{proposition}
The proposition remains valid if ``half-Eulerian'' is replaced by 
$k$-Eulerian, and each $f_S$ is replaced by $(2k)^{n-|S|}f_S$.

\vspace*{6pt}

\begin{proof}
We first show the inequality is not a convolution of lower rank inequalities.
In $L^k$-vector form the convolution satisfies the rule
$L^{k,i+1}_T * L^{k,j+1}_V = 2kL^{k,i+j+2}_{T\cup(V+i+1)}$ (see
equation~\ref{convol}).
So the convolution of linear expressions for ranks $i+1$ and $j+1$ with
$n=i+j+1$ gives a linear combination of $L^{k,n+1}_S$ involving only
subsets $S\subseteq[1,n]$ not containing $i+1$.
Since each element of $[1,7]$ occurs in some set $S$ in the 
inequality~\ref{rank8ineq}, it is not a convolution of lower rank
inequalities.

We now show that the inequality determines a facet of the cone.
Billera and Hetyei list the facet inequalities for the general graded
cone up through rank 5.
The inequality of the proposition comes from applying one of the
rank 5 Billera-Hetyei inequalities to the rank-selected subposet 
$P_{\{1,3,5,6\}}$ of an arbitrary half-Eulerian poset $P$.
To check it is a facet of the half-Eulerian cone, we give twenty linearly
independent limiting normalized $L^{1/2}$-vectors of half-Eulerian posets, for 
which the inequality holds with equality.
The first sixteen posets are Billera-Hetyei limit posets determined by
interval systems as in the following table.

\vspace*{6pt}

\begin{tabular}{||l|c||l|c||l|c||} \hline\hline
$P_1$ & $\emptyset$ & $P_7$ & $[1,2][5,6]$ & $P_{12}$ & $[1,4][6,7]$ \\ \hline
$P_2$ & $[1,2]$ & $P_8$ & $[1,2][3,4][5,6]$ & $P_{13}$ & $[4,5][6,7]$ \\ \hline 
$P_3$ & $[2,3]$ & $P_9$ & $[3,6]$ & $P_{14}$ & $[2,3][4,5][6,7]$ \\ \hline
$P_4$ & $[3,4]$ & $P_{10}$ & $[6,7]$ & $P_{15}$ & $[1,2][4,7]$ \\ \hline
$P_5$ & $[1,2][3,4]$ & $P_{11}$ & $[1,2][6,7]$ & $P_{16}$ & $[2,7]$ \\ \hline
$P_6$ & $[2,3][4,5]$ & \multicolumn{4}{c||}{} 
\\ \hline\hline
\end{tabular}

\vspace*{6pt}

The next three limit posets are obtained from the rank 7 Extremes~2, 3 and~4 of
\cite[Theorem 4.8]{Bayer-Hetyei} by inserting a single new element of rank 1,
shifting the old elements up one rank.

To describe the last sequence of posets, let us (re)introduce the
following generalization of the operator $D^r$. Given a graded poset $P$
of rank \mbox{$n+1$} denote by $D^r_{[u,v]}(P)$ the poset
obtained from $P$ by replacing each $x\in P$ satisfying $\rho(x)\in [u,v]$ 
with $r$ elements $x_1,x_2,\ldots, x_r$ (keep every $y\in P$
satisfying $\rho (y)\not\in [u,v]$ unchanged), and by setting the
following order relations. The $([1,n]\setminus [u,v])$-rank-selected
subposet of $P$ and of $D^r_{[u,v]}(P)$ are identical. For $x,y\in P$
satisfying $\rho(x)\in [u,v]$ and $\rho (y)\not \in [u,v]$ set 
$x_i<y$ or $x_i>y$ in $D^r_{[u,v]}(P)$ if and only of the same
relation holds between $x$ and $y$ in $P$. Finally for $x,y\in P$
satisfying $u\leq \rho (x)<\rho(y)\leq v$ set $x_i<y_j$
in $D^r_{[u,v]}(P)$ if and only if $i=j$ and $x<y$ in $P$.

For example, Figure \ref{F_example} shows $D^2_{[1,2]}(C_4)$ where $C_4$ is
a chain of rank $4$. Note that for a graded poset $P$ of rank $n+1$ the graded
poset $D^r (P)$ is isomorphic to $D^r_{[1,1]}D^r_{[2,2]}\ldots
D^r_{[n,n]}(P)$. The same notation is used in \cite{Bayer-Hetyei}.

\begin{figure}[h]
\setlength{\unitlength}{0.00063300in}%
\begingroup\makeatletter\ifx\SetFigFont\undefined
\def\x#1#2#3#4#5#6#7\relax{\def\x{#1#2#3#4#5#6}}%
\expandafter\x\fmtname xxxxxx\relax \def\y{splain}%
\ifx\x\y   
\gdef\SetFigFont#1#2#3{%
  \ifnum #1<17\tiny\else \ifnum #1<20\small\else
  \ifnum #1<24\normalsize\else \ifnum #1<29\large\else
  \ifnum #1<34\Large\else \ifnum #1<41\LARGE\else
     \huge\fi\fi\fi\fi\fi\fi
  \csname #3\endcsname}%
\else
\gdef\SetFigFont#1#2#3{\begingroup
  \count@#1\relax \ifnum 25<\count@\count@25\fi
  \def\x{\endgroup\@setsize\SetFigFont{#2pt}}%
  \expandafter\x
    \csname \romannumeral\the\count@ pt\expandafter\endcsname
    \csname @\romannumeral\the\count@ pt\endcsname
  \csname #3\endcsname}%
\fi
\fi\endgroup
\begin{center}
\begin{picture}(1066,2602)(818,-5846)
\thicklines
\put(901,-4523){\circle*{150}}
\put(901,-4973){\circle*{150}}
\put(1351,-5423){\circle*{150}}
\put(1351,-4073){\circle*{150}}
\put(1351,-3623){\circle*{150}}
\put(1801,-4973){\circle*{150}}
\put(1801,-4523){\circle*{150}}
\put(1300,-5798){\makebox(0,0)[lb]{\smash{\SetFigFont{12}{14.4}{rm}$\hat{0}$}}}
\put(1300,-3436){\makebox(0,0)[lb]{\smash{\SetFigFont{12}{14.4}{rm}$\hat{1}$}}}
\put(1351,-3661){\line( 0,-1){450}}
\put(1351,-4111){\line(-1,-1){450}}
\put(901,-4561){\line( 0,-1){450}}
\put(901,-5011){\line( 1,-1){450}}
\put(1351,-5461){\line( 1, 1){450}}
\put(1801,-5011){\line( 0, 1){450}}
\put(1801,-4561){\line(-1, 1){450}}
\end{picture}
\end{center}
\caption{$D^2_{[1,2]}(C_4)$}
\label{F_example}
\end{figure}

Let $N$ be an arbitrary positive integer, and $C_8$ be a chain of rank $8$.
Consider now the following four graded posets.
\begin{eqnarray*}
P^I(N) &=& D^{N+1}_{[1,2]}D^{N+1}_{[2,3]}D^{N+1}_{[4,5]}D^{N}_{[1,7]}
(C_8)\\
P^{II}(N) &=& D^{N^2}_{[1,3]}D^{N+1}_{[1,5]}D^{N}_{[1,7]} (C_8)\\
P^{III}(N) &=& D^{N^2-N+2}_{[1,4]}D^{N+2}_{[4,5]}D^{N}_{[6,7]} (C_8)\\
P^{IV}(N) &=& D^{N+2}_{[1,2]}D^{N^3-N^2+2}_{[2,7]}(C_8)\\
\end{eqnarray*}
The $\{4,5,6,7\}$-rank-selected subposets of $P^I(N)$ and $P^{II}(N)$ are
both isomorphic to $D^{N+1}_{[1,2]}D^{N}_{[1,4]}(C_5)$, where $C_5$ is a
chain of rank $5$; the $\{6,7\}$-rank-selected subposets of $P^I(N)$,
$P^{II}(N)$, and $P^{III}(N)$ are all isomorphic to 
$D^{N}_{[1,2]}(C_3)$ where $C_3$ is a chain of rank $3$. Let $P(N)$ be
the graded poset of rank $8$ obtained from $P^I(N)$,
$P^{II}(N)$, $P^{III}(N)$, and $P^{IV}(N)$ by performing the following
identifications:

-identify the bottom element $\hat{0}$ of all four posets,

-identify the top element $\hat{1}$ of all four posets,

-identify $P^I(N)_{\{4,5,6,7\}}$ with
$P^{II}(N)_{\{4,5,6,7\}}$,

-identify $P^I(N)_{\{6,7\}}$ with
$P^{III}(N)_{\{6,7\}}$.

Figure \ref{F_4d1} indicates how the four posets are identified,
in a schematic way.

\begin{figure}[h]
\setlength{\unitlength}{0.00063300in}%
\begingroup\makeatletter\ifx\SetFigFont\undefined
\def\x#1#2#3#4#5#6#7\relax{\def\x{#1#2#3#4#5#6}}%
\expandafter\x\fmtname xxxxxx\relax \def\y{splain}%
\ifx\x\y   
\gdef\SetFigFont#1#2#3{%
  \ifnum #1<17\tiny\else \ifnum #1<20\small\else
  \ifnum #1<24\normalsize\else \ifnum #1<29\large\else
  \ifnum #1<34\Large\else \ifnum #1<41\LARGE\else
     \huge\fi\fi\fi\fi\fi\fi
  \csname #3\endcsname}%
\else
\gdef\SetFigFont#1#2#3{\begingroup
  \count@#1\relax \ifnum 25<\count@\count@25\fi
  \def\x{\endgroup\@setsize\SetFigFont{#2pt}}%
  \expandafter\x
    \csname \romannumeral\the\count@ pt\expandafter\endcsname
    \csname @\romannumeral\the\count@ pt\endcsname
  \csname #3\endcsname}%
\fi
\fi\endgroup
\begin{center}
\begin{picture}(4921,5077)(1913,-6408)
\thicklines
\put(4051,-1823){\circle*{150}}
\put(2701,-2723){\circle*{150}}
\put(2701,-3173){\circle*{150}}
\put(4051,-4523){\circle*{150}}
\put(4051,-4973){\circle*{150}}
\put(4051,-5423){\circle*{150}}
\put(2701,-4523){\circle*{150}}
\put(2701,-4073){\circle*{150}}
\put(2701,-3623){\circle*{150}}
\put(4051,-5873){\circle*{150}}
\put(2701,-5423){\circle*{150}}
\put(2701,-4973){\circle*{150}}
\put(6751,-4073){\circle*{150}}
\put(6751,-4523){\circle*{150}}
\put(6751,-4973){\circle*{150}}
\put(6751,-2723){\circle*{150}}
\put(6751,-3173){\circle*{150}}
\put(6751,-3735){\circle*{150}}
\put(5401,-4523){\circle*{150}}
\put(5401,-4973){\circle*{150}}
\put(5401,-5423){\circle*{150}}
\put(6751,-5423){\circle*{150}}
\put(5401,-3623){\circle*{150}}
\put(5401,-4073){\circle*{150}}
\put(3153,-4898){\makebox(0,0)[lb]{\smash{\SetFigFont{12}{12}{rm}$P^{II}(N)$}}}
\put(5751,-4898){\makebox(0,0)[lb]{\smash{\SetFigFont{12}{12}{rm}$P^{IV}(N)$}}}
\put(1913,-4898){\makebox(0,0)[lb]{\smash{\SetFigFont{12}{12}{rm}$P^I(N)$}}}
\put(4000,-1523){\makebox(0,0)[lb]{\smash{\SetFigFont{12}{12}{rm}$\hat{1}$}}}
\put(4000,-6360){\makebox(0,0)[lb]{\smash{\SetFigFont{12}{12}{rm}$\hat{0}$}}}
\put(4051,-5911){\line( 3, 1){1350}}
\put(5401,-5461){\line( 0, 1){1800}}
\put(5401,-3661){\line(-6, 1){2700}}
\put(4051,-5911){\line( 6, 1){2700}}
\put(6751,-5461){\line( 0, 1){2700}}
\put(6751,-2761){\line(-3, 1){2700}}
\put(4051,-5911){\line( 0, 1){1350}}
\put(4051,-4561){\line(-3, 1){1350}}
\put(4401,-4898){\makebox(0,0)[lb]{\smash{\SetFigFont{12}{12}{rm}$P^{III}(N)$}}}
\put(4051,-5911){\line(-3, 1){1350}}
\put(2701,-5461){\line( 0, 1){2813}}
\put(2701,-2648){\line( 5, 3){1339.853}}
\end{picture}
\end{center}
\caption{$P(N)$}
\label{F_4d1}
\end{figure}

Straightforward calculation shows that $P(N)$ is a half-Eulerian poset,
for each positive $N$. 
Furthermore the normalized $L^{1/2}$-vectors, 
$(L^{1/2,8}_S(P(N))/N^4)$, converge.

The rows of the matrix below are the normalized $L^{1/2}$-vectors of the twenty
limit posets.
In the columns are the values of $L^{1/2,8}_S$ (divided by the appropriate
power of $N$), with the sets $S$ in the order
$\emptyset$, $\{1,2\}$, $\{2,3\}$, $\{3,4\}$, $\{1,2,3,4\}$, $\{4,5\}$, 
$\{1,2,4,5\}$, $\{2,3,4,5\}$, $\{5,6\}$, $\{1,2,5,6\}$, $\{2,3,5,6\}$, 
$\{3,4,5,6\}$, $\{1,2,3,4,5,6\}$, $\{6,7\}$, $\{1,2,6,7\}$, $\{2,3,6,7\}$, 
$\{3,4,6,7\}$, $\{1,2,3,4,6,7\}$, $\{4,5,6,7\}$, $\{1,2,4,5,6,7\}$, 
$\{2,3,4,5,6,7\}$. 
It is easy to check the rows are linearly independent.

\scriptsize{
$$ \left[
\begin{array}{rrrrrrrrrrrrrrrrrrrrr}
 1 & 0 & 0 & 0 & 0 & 0 & 0 & 0 & 0 & 0 & 0 & 0 & 0 & 0 & 0 & 0 & 0 & 0 & 0 & 0 & 0\\
 1 & -1 & 0 & 0 & 0 & 0 & 0 & 0 & 0 & 0 & 0 & 0 & 0 & 0 & 0 & 0 & 0 & 0 & 0 & 0 & 0\\
 1 & 0 & -1 & 0 & 0 & 0 & 0 & 0 & 0 & 0 & 0 & 0 & 0 & 0 & 0 & 0 & 0 & 0 & 0 & 0 & 0\\
 1 & 0 & 0 & -1 & 0 & 0 & 0 & 0 & 0 & 0 & 0 & 0 & 0 & 0 & 0 & 0 & 0 & 0 & 0 & 0 & 0\\
 1 & -1 & 0 & -1 & 1 & 0 & 0 & 0 & 0 & 0 & 0 & 0 & 0 & 0 & 0 & 0 & 0 & 0 & 0 & 0 & 0\\
 1 & 0 & -1 & 0 & 0 & -1 & 0 & 1 & 0 & 0 & 0 & 0 & 0 & 0 & 0 & 0 & 0 & 0 & 0 & 0 & 0\\
 1 & -1 & 0 & 0 & 0 & 0 & 0 & 0 & -1 & 1 & 0 & 0 & 0 & 0 & 0 & 0 & 0 & 0 & 0 & 0 & 0\\
 1 & 0 & 0 & 0 & 0 & 0 & 0 & 0 & 0 & 0 & 0 & -1 & 0 & 0 & 0 & 0 & 0 & 0 & 0 & 0 & 0\\
 1 & -1 & 0 & -1 & 1 & 0 & 0 & 0 & -1 & 1 & 0 & 1 & -1 & 0 & 0 & 0 & 0 & 0 & 0 & 0 & 0\\
 1 & 0 & 0 & 0 & 0 & 0 & 0 & 0 & 0 & 0 & 0 & 0 & 0 & -1 & 0 & 0 & 0 & 0 & 0 & 0 & 0\\
 1 & -1 & 0 & 0 & 0 & 0 & 0 & 0 & 0 & 0 & 0 & 0 & 0 & -1 & 1 & 0 & 0 & 0 & 0 & 0 & 0\\
 1 & 0 & 0 & 0 & -1 & 0 & 0 & 0 & 0 & 0 & 0 & 0 & 0 & -1 & 0 & 0 & 0 & 1 & 0 & 0 & 0\\
 1 & 0 & 0 & 0 & 0 & -1 & 0 & 0 & 0 & 0 & 0 & 0 & 0 & -1 & 0 & 0 & 0 & 0 & 1 & 0 & 0\\
 3 & 0 & -1 & -1 & 0 & -1 & 0 & 1 & -1 & 0 & 0 & 0 & 0 & -1 & 0 & 0 & 0 & 0 & 1 & 0 & 0\\
 3 & 0 & -2 & 0 & 0 & -1 & 0 & 1 & -1 & 0 & 1 & -1 & 0 & -1 & 0 & 0 & 0 & 0 & 1 & 0 & 0   \\
 3 & 0 & -1 & -1 & 0 & -1 & 0 & 1 & 0 & 0 & 0 & -1 & 0 & -2 & 0 & 0 & 1 & 0 & 1 & 0 & 0  \\
 1 & -1 & 0 & 0 & 0 & 0 & 0 & 0 & 0 & 0 & 0 & 0 & 0 & 0 & 0 & 0 & 0 & 0 & -1 & 1 & 0\\
 1 & 0 & -1 & 0 & 0 & -1 & 0 & 1 & 0 & 0 & 0 & 0 & 0 & -1 & 0 & 1 & 0 & 0 & 1 & 0 & -1\\
 1 & 0 & 0 & 0 & 0 & 0 & 0 & 0 & 0 & 0 & 0 & 0 & 0 & 0 & 0 & 0 & 0 & 0 & 0 & 0 & -1\\
 4 & -2 & -1 & 0 & -1 & -2 & 1 & 1 & 0 & 0 & 0 & 0 & 0 & -1 & 0 & 0 & 0 & 1 & 1 & 0 & -1
\end{array}
\right] $$}

\end{proof}


\end{document}